\documentclass[ptm,10pt,reqno]{amsart} 
\usepackage{amsmath,amssymb,rawfonts}
\usepackage{multicol}
\usepackage{mathptmx}
\usepackage{tikz}
\usepackage{pgfplots}
\usepackage{color}
\usepackage{hyperref}
\usepackage[utf8]{inputenc}
\usepackage[all]{xy}
\usetikzlibrary{matrix,arrows}

\newtheorem{theorem}{Theorem}[section]
\newtheorem{lemma}[theorem]{Lemma}
\newtheorem{defi}[theorem]{Definition}

\newtheorem{aplication}[theorem]{Application}

\newtheorem{remark}[theorem]{Remark}

\begin{document} 
\title{A context for  Paris-Harrington combinatorial Principle}
\author{Joel Torres Del valle}
\address{Universidad de Antioquia, FCEyN, Instituto de Matemáticas.}
\email{jtorresdv1@gmail.com}
\maketitle

\noindent For a long time there was the pretence of conceiving mathematics as an idealization of the palpable world, and proceeding on it as one would between objects of the real world. At the end of the 19th century paradoxes began to appear in the young Set Theory of Georg Cantor and thus mathematics, which bore the unworthy title of an exact science, began to fade away. It was observed, then, that the beautiful building stood on quicksand and staggered to the sound of the lightest wind.

Then began a program of foundation of mathematics, which was intended to build solid foundations on which to stop the mathematical building. Of the various currents of mathematical philosophy that addressed the problem, we can point to Hilbert's Formalism.  The formalist program was mainly driven by David Hilbert, and his final requirement was proof of the non-contradiction of mathematics. Of course, we take for granted, the dream of the completeness of the formal systems upon which mathematics would be based. For David Hilbert, logic and mathematics are theories of form and not of meaning, \cite{panorama}. That is to say, all the tests were to be carried out by means of fixed rules on the handling of symbols, without taking into account at any time the meaning of the symbols. Thought that is clear when this says: mathematics is a game with very simple rules, which leave meaningless marks on paper. 

The zeal for a proof of non-contradiction comes after the works of Gottlob Frege (1848-1925) \cite{frege} could be deduced  Russell's paradox, which in Peano's symbolism could be expressed as\footnote{Take from Bertrand Russell - Gottlob Frege correspondence.} $w=\mathrm{cls}\cap x\backepsilon(x\sim \epsilon x).\supset:w\in w.=.w\sim\epsilon w$, which, in Russell's own words, corresponds to: let $w$ be the predicate: to be a predicate that cannot be predicated of itself. Can $w$ be a predicated of itself? From each answer its opposite follows. Therefore we must conclude that $w$ is not a predicate. Likewise there is no class (as a totality) of those classes which, each taken as a totality, do not belong to themselves. 

And I cannot let pass, of course, the height with which Frege assumes the terrifying news that his work allowed the appearance of paradoxes and, full of scientific humility, he adds a note to the second volume of his work \cite{frege-ii} which was already in print, in which he comments on Russell's discovery in the first volume of his work. 

In the year 1910 appears the first volume of {\bf Principia Mathematica}, in which Bertrand Russell (1972-1970) and Alfred North Whitehead (1861-1947), build much of mathematics without paradoxes or apparent contradictions. The question then remained as to whether the system was complete, given that everything seemed to indicate consistency. The first steps towards a clarification of this requirement were taken in a direction that seemed to discern an light at the end of the tunnel. Firstly in 1930, in his doctoral thesis  \cite{godel} Kurt Gödel (1906-1978) proved the completeness of First-Order Logical Calculus. However, by these paradoxes of luck, it is the same Gödel who a year later, in 1931 proved that if we add to Peano's Axioms for arithmetic  all the logic of Principia Mathematica, we do not obtain a system from which all the truths about natural numbers can be demonstrated,  \cite{godelii}. Specifically (see {\it Some mathematical results on completeness and consistency}, Kurt Gödel, 1931)

\begin{itemize}
\item[1.] The $S$ system is not complete, that is, in it there are sentences $\varphi$, such that neither $\varphi$ nor $\neg\neg\varphi$ are probable  and, in particular, there are undecidable problems with the simple structure $\exists xFx$, where $x$ varies over natural numbers and $F$ is a property of natural numbers.
\item[2.] Even if we admit all the logical means of Principia Mathematica in metamathematics there is no proof of consistency for $S$. Therefore, a consistency test for the $S$ system can only be carried out with the help of inference modes that are not formalized in the $S$ system.
\item[3.] Not even adding to $S$ a finite amount of axioms in such a way that the extended system remains $\omega$-consistent.
\end{itemize}

However, these phrases initially were not something that a mathematician would ask himself, for example, a mathematician would ask himself is Goldbach's conjecture true or false? etc. However, through a process today known as Gödel Coding, Gödel constructs an explicit phrase, from the Theory of Numbers (but which does not arise naturally in it) that cannot be proved or refuted. Something like a mathematical formulation of the liar paradox\footnote{Gödel demonstrates that natural numbers are strong enough to codify all the truths of the $S$ system and about the S.}.   So the question remains, how can we make these phrases more mathematical? There appeared Paris-Harrington result. 

In previous years had already established such results in Euclidean Geometry (in the 19th century) and Set Theory, namely the V Postulate of Euclid by János Bolyai and Nikolái Lobachesky  and the Continuum Hypothesis, by Kurt Gödel and Paul Cohen. On the V Postulate of Euclides, it was proved that it is independent of the Euclid's axiomatic system, and that in the same way geometries can be constructed where it is valid and where it is not, which gives the appearance to non-Euclidean Geometries. The same thing happened with the continuum of Cantor. 

In this pages we present a context for Paris-Harrington Principle. We follow \cite{maker} for Section 3, \cite{kaye} for Section 4; for Section 1 we follow \cite{luis-silva}.  For some background in Combinatorial Set Theory we use \cite{keku}, \cite{luis-silva}. It is important to make clear, that in the formulation of PH which we study, we just consider partitions of sets like $m:=\{0, 1, ..., m-1\}$. A general formulation can be done following the same argument that we present here for $H\in\mathrm{Par}([m]^{n}, c)$.

\section{The Paris Harrington principle}

\subsection{The König's Lemma}\hspace{0.02cm}\\

 We start with some background in Combinatorial Set Theory. For the proof of Infinite Ramsey Theorem (and so for that of PH) we need the König's Lemma which we prove here. For this section we follow \cite{luis-silva, keku, maker}.

\begin{defi}
\rm A \textbf{tree} is a partially ordered set  $(T, <_{T})$ such that for all $t\in T$, the set  $\widehat{t}:=\{s\in T:s<_{T}t\}$ is well ordered. A  \textbf{branch} of a tree $T$ is a co-final chain (a  linearly ordered set) maximal of $T$. A \textbf{trajectory} of $T$ is a chain of $T$ which is also a initial segment of $T$.
\end{defi}
 
\begin{multicols}{3}

\[
\xymatrix{ &\cdot & &\cdot \\
     \cdot & & \cdot\ar@{-}[ur]\ar@{-}[ul] & \\
           & \cdot\ar@{-}[ur]\ar@{-}[ul] & & \\}
\]

\begin{center}
\textbf{Fig. 1.} 
\end{center}

\[
\xymatrix{ & & &\cdot \\
      & & \cdot\ar@{-}[ur]  & \\
           & \cdot\ar@{-}[ur]& & \\}
\]

\begin{center}
\textbf{Fig. 2.} 
\end{center}

\[
\xymatrix{ &\cdot  & &\cdot \\
     & & \cdot\ar@{-}[ur]\ar@{-}[ul]  & \\
           & \cdot\ar@{-}[ur]& & \\}
\]

\begin{center}
\textbf{Fig. 3} 
\end{center}
\end{multicols}

\noindent In Figure 1. it is drawn a tree, in Figure 2. a trajectory and in Figure 3. a branch of a tree.  

\begin{defi}
\rm
A \textbf{finite branching tree } is a partially ordered set  $(T, <_{T})$ such that: there exists  $r\in T$ such that  $r<_{T} x$ for all $x\in T$. If $x\in T$, then $\{y:y<_{T}x\}$ is finite and linearly ordered by $<_{T}$ and finally if $x\in T$, each finite set (maybe empty) $\{y_{1}, ..., y_{n}\}$ of incomparable elements such that each $y_{i}>x$ and each $z>x$, then $z>y_{i}$ for some $i$.
\end{defi} 

It is known that  any well ordered set is isomorphic to some unique ordinal.  This ordinal is known as \textbf{type ordinal} of the set. The  \textbf{height} $\mathrm{hgt}(t)$ of $t$ in $T$ is the ordinal type of  $\widehat{t}$. The \textbf{level} $\alpha$ of $T$ is the set $T_{\alpha}:=\{t\in T:\mathrm{hgt}(t)=\alpha\}$. The \textbf{height} of  $T$ is $\min\{\alpha:T_{\alpha}=\emptyset\}$.

\begin{defi}
\rm
Let be $\theta$ an ordinal and $\lambda$ a cardinal. A tree $T$ is a  $(\theta, \lambda)$-\textbf{tree} if: $(\forall\alpha<\theta)(T_{\alpha}\neq\emptyset)$, $T_{\theta}=\emptyset$ and $(\forall\alpha<\theta)(\mathrm{card}(T_{\alpha})<\lambda)$. An $\aleph_{0}$-\textbf{tree}, is simple a  $(\aleph_{0}, \aleph_{0})$-tree. 
\end{defi}

\begin{lemma}[Köning's Lemma] Any $\aleph_{0}$-tree has a branch which intersects any levels.
\end{lemma} 

In fact, let  $T$ be an $\aleph_{0}$-tree. By induction on  $n<\omega$, let us choice  $t_{n}\in T_{n}$ such that $T^{t_{n}}$ is infinite and $t_{n}<_{T}t_{n+1}$. Then, $\{t_{n}:n<\omega\}$ is a co-final branch of $T$.

\subsection{Partitions theorems}\hspace{0.02cm}\\

The next step is to prove Infinite and Finite Ramsey theorems for partitions. 

\begin{defi}
\rm
 Let $\sigma$ be a cardinal, $[I]^{k}$ the set of those subset of $I$ of size $k$. A map $P:[I]^{k}\to\sigma$ is called a \textbf{partition} of $[I]^{k}$ in $\sigma$ pieces. A set $Y\subset I$ is called  \textbf{homogeneous for the partition} $P$, if $P$ is constant on $[H]^{k}$.
 \end{defi}

We note $\kappa\to(\lambda)_{\sigma}^{n}$ to indicate that: whenever $P:[\kappa]^{n}\to\sigma$ is a partition of $[\kappa]^{n}$ into $\sigma$ pieces, then there exists $Y\subset\kappa$ of size at most $\lambda$ homogeneous for $P$. If it also $\mathrm{card}(Y)\geq\min(Y)$ we call $Y$  \textbf{relatively large}, and we note $\xymatrix{
\kappa\ar[r]_{*}&(\lambda)_{\sigma}^{n}\\}$. 

\begin{theorem}[Infinite Ramsey Theorem]
For any natural numbers $n, k$ it holds that $\aleph_{0}\to(\aleph_{0})_{k}^{n}$.
\end{theorem} 

We proceed by induction on  $n$. For  $n=0$ there it nothing to prove because  $f$ is constant on $[A]^{0}=\{\emptyset\}$. Let $n>0$.

 We define a decreasing sequence recursively  $A_{0}\supset A_{1}\supset ...$ of infinite sets of  $A$ and a sequence  $a_{0}, a_{1}, ...$ of elements of  $A$ with $a_{i}\in A_{j}$ only if  $i\geq j$.

 We start with $A_{0}=A$. Let us suppose that  $A_{i}$ is given. Let it be $a_{i}\in A_{i}$ arbitrary. We define $f_{i}:[A_{i}\setminus\{a_{i}\}]^{n-1}\to m$ by $f_{i}(b)=f(\{a_{i}\}\cup b)$. As  $A_{i+1}$, we choose $f_{i}$-homogeneous of $A_{i}\setminus\{a_{i}\}$.

Let it be $m_{i}$ the value that takes  $f_{i}$ in $[A_{i+1}]^{n-1}$. Then, for each $k<m$ the set $B=\{a_{i}:m_{i}=k\}$ is $f$-homogeneous: each set of  $n$ elements $c$ of $B$ is of the form $\{a_{i}\}\cup b$ for some $b\in [A_{i+1}]^{n-1}$. We have $f(c)=f_{i}(b)=k$. Then there exists  $k<m$ such that $m_{i}=k$ for  a infinite quantity of $i\in\omega$. $B$ is infinite for this $k$. This completes the proof. 

\begin{aplication}[Cited from \cite{maker}]
\rm
Any sequence $\{r_{i}\}_{i\in\mathbb{N}}$ of real numbers has a monotonic subsequence.  In fact, let be  $f:[\mathbb{N}]^{2}\to 3$ defined by

\[
f(\{i, j\}):=\left\{\begin{array}{cccc}
0&i<j&\mathrm{and}&r_{i}<r_{j}\\
1&i<j&\mathrm{and}&r_{i}=r_{j}\\
2&i<j&\mathrm{and}&r_{i}>r_{j}\\
\end{array}\right.
\]

\noindent By the Infinite Ramsey Theorem, there is $Y\subset \mathbb{N}$ with $\mathrm{card}(Y)=\aleph_{0}$ such that $f$ is constant on $[Y]^{2}$.  Let $j_{0}<j_{1}<...$ list $Y$. There is $c\in 3:=\{0, 1, 2\}$ such that $f(\{i_{m}, j_{n}\})=c$ for $m<n$. If $c=0$, the sequence $\{r_{j_{k}}\}_{k\in Y}$ is increasing, if $c=1$ it is constant, and $c=2$ it is decreasing. 
\end{aplication}

\begin{theorem}[Finite Ramsey Theorem]
For any natural numbers $n, r, k$ there exists a natural number $l$ such that $l\to(r)_{k}^{n}$.
\end{theorem} 

In fact, suppose there is not such $l$. For each $l<\omega$, let  be  $T_{l}$ the set defined by:  
\begin{small}
$T_{l}:=\left\lbrace f:[l]^{n}\to k:\right.$ there not exists a set of  $l$ of size $r$ homogeneous for $\left.f\right\rbrace$.
\end{small}

 Clearly,   $T_{l}$ is finite, since there are finite partitions of finite sets. Let be $f\in T_{l+1}$, then there is an unique $g\in T_{l}$ such that $g\subset f$. Hence, if we order the set  $
T=\bigcup_{l<\omega}T_{l}$  by inclusion, we obtain a finite tree. Each $T_{j}\neq\emptyset$. Then, we obtain a tree of finite branching. By Köning's Lemma, we can find  $f_{0}\subset f_{1}\subset ...$  with each $f_{i}\in T_{i}$.  Let be $f=\bigcup f_{i}$, then $f:[\mathbb{N}]^{n}\to k$. For Infinite Ramsey Theorem, there is  a $X\subset \mathbb{N}$ infinite homogeneous for $f$. Let be $x_{1}, ..., x_{r}$, the first $r$ elements of $X$ and let be $s>x_{r}$, then $\{x_{1}, ..., x_{r}\}$ is homogeneous for $f_{s}$. This is a contradiction.

\begin{aplication}
\rm $6\to (3)_{2}^{2}$. That is, for any partition $f:[6]^{2}\to 2$ there exists a subset of 6 with  size 3 such that $f:[3]^{2}\to 2$ is constant. As a consequence of this, we have that in a meeting of six people, three o them know each other, or three of them don't know each other. For this, consider the set $p=\{p_{1}, p_{2}
, p_{3}, p_{4}, p_{5}, p_{6}\}$ of people and $f$ the function 

\[
\begin{array}{cccc}
f:&[p]^{2}&\hspace{0.66cm}\longrightarrow\hspace{0.5cm}2& \\
  &\{p_{i}, p_{j}\}&\longmapsto&\left\{\begin{array}{ccc}
  1&\mathrm{if}&p_{i}\hspace{0.2cm}\mathrm{knows}\hspace{0.2cm}p_{j}\\
  0&\mathrm{if}&p_{i}\hspace{0.2cm}\mathrm{doesn't}\hspace{0.2cm}\mathrm{know}\hspace{0.2cm}p_{j}
  \end{array}\right.
\end{array}
\] 

\noindent By Finite Ramsey theorem, there exists $p'\subset p$ with $\mathrm{card}(p')=3$ such that $f:[p']^{2}\to 2$ is constant. Then $f([p']^{2})=0$ or $f([p']^{2})=1$. So, we obtain the result.  
\end{aplication}

\begin{theorem}[Paris-Harrington Principle] 
For any natural numbers $n, r, k$ there exists a natural number $m$ such that $\xymatrix{m\ar[r]_{*}&(k)_{r}^{n}\\}$.
\end{theorem} 

Paris-Harrington result was the first mathematical independence in arithmetic result since Gödel's Theorems. And it was the answer to the question: How can we make independence results mathematical? Because of Gödel presented an explicit sentence which is independent of PA, but these didn't appear naturally from the number theory, nor it was a metamathematical statement which asserted it's own independence;  it was important to answer the question above.  

\begin{proof} Let be  $n, r$ and  $k$  natural numbers. Suppose that there not exists such $m$. Let  $P$ be a counterexample for $m$. If $P$ is a partition of $[m]^{n}$ in  $r$ parts with any homogeneous subset of size at most $k$. We can see the set of counterexamples as a infinite tree of finite branching, that is to say, if $P$ and $P'$ are counterexamples for  $m$ and $m'$ respectively, we put $P$ below $P'$ in our  tree only if $m<m'$ and  $P$ is a restriction of  $P'$ to $[m]^{n}$.

By  König's Lemma there is a $P:[\omega]^{n}\to r$ such that for any $m$, the restriction of $P$ to $[m]^{n}$ is a counterexample for $m$. By Infinite Ramsey Theorem, there exists a $H\subset \omega$ infinite homogeneous for $P$. But then, by taking $m$ enough large (comparing with $k$ and $\min H$) we see that  $H\cap m$ is, after all, an homogeneous set relatively large for $P\upharpoonright [m]^{n}$ of size at least $k$.
\end{proof}

Finite Ramsey Theorem is an assertion about integers, so we would prefer to give a proof for it without using infinite methods as we do here. Our approach is given like in \cite{maker, luis-silva}, but in other works a finite approach is given. On the other hand since PH  looks like Finite Ramsey Theorem we would expect the same occurs with it, but this is not the case: we can't carry out PH proof without inifinite methods, so we can't neither formalize the proof of PH in the language of arithmetic as we can do with that of Finite Ramsey Theorem.

By using the argument of the Second Incompleteness
Theorem, on 70's years, Jeff Paris and Leo Harrington, showed
that an strength of Finite Ramsey Theorem is not provable from Peano
Arithmetic, PA.  In \cite{paris} Paris and Harrington introduce a certain theory 
$T$, and they show that $\mathrm{Con}(T)\to \mathrm{Con}(\mathrm{PA})$ is theorem of PA, then the proof of the independence
of PH is complete when they show that $\mathrm{PH}\to \mathrm{Con}(T)$ is also a theorem of
PA, since for the Second Incompleteness Theorem it would impossible that $\mathrm{PA}\vdash \mathrm{PH}$.

\section{Encoding finite subsets and partitions}\hspace{0.02cm}\\

In 1931, Kurt Gödel, introduced a process which nowadays is known as Gödel numbering, and it permitted him to translate into the language $\mathcal{L}_{A}:=\{0, 1, +, \cdot, <\}$ of Arithmetic some interesting and important metamathematical properties on Number Theory, such as \emph{Consistency}, etc. By using this, Gödel proved his celebrated Incompleteness Theorems; these ones, represented the end of the Formalist Dream. Gödel numbering also permits us to express any properties of Finite Combinatorics and Partitions, into the language $\mathcal{L}_{A}$.

\subsection{Gödel codes for finite sets.}\hspace{0.02cm}\\

For this section I mainly follow \cite{kaye}. Identify each natural number with their set of predecessors, i.e., $m:=\{0, 1, ..., m-1\}$.  We will work with sets written in  increasing order. The \textbf{Gödel code} for the set $a=\{a_{0}, a_{1}, ..., a_{n}\}$ is defined by $:=2^{a_{0}}3^{a_{1}}\cdot\cdot\cdot p_{n}^{a_{n}}=\prod_{i=0}^{n}p_{i}^{a_{i}}$. Here, $p_{n}$ represents the $n$-prime number. For the uniqueness of prime factorization, the Gödel code of a set $a$ is a unique number. It is clear that the number of prime factors in the code of $a$ is the same number of elements in $a$. 

\begin{defi}\rm  Let us write $\mathrm{Prime}(p)$ to mean  `$p$ is a prime number'. This is given the formula:  

\[
\mathrm{Prime}(p):=p\neq 0\wedge p\neq 1\wedge\forall x<p![x|p\to p=x\vee x=1], x|p:=\exists z[xz=p].
\]

\end{defi}

\begin{defi}\rm 
Let $g(x, y)$ be a function, we define the $\mu$-\textbf{bounded operator} by: $f(x, x)=\mu y<x[g(x, y)=0]$ to be: the least $y<x$ such that $g(x, y)=0$ if there exists such $y$, and $0$ otherwise.  `$p_{n}$ is the $n$-th prime number' is defined by  $p_{n}:=p_{0}=2, p_{n+1}:=\mu x<p_{n}!+1[p_{n}<x\wedge\mathrm{Prime}(x)]$.
\end{defi}

\begin{defi}\rm We define `Seq' to be the set of those numbers which codifies sequences as it follows:

\[
a\in \mathrm{Seq}\leftrightarrow a=1\vee (a>1\wedge\forall x\leq a[p_{x+1}|a\to p_{x}|a]).
\]

Let be $a\in $Seq, we define `Lgh($a$)' to be the \textbf{length} of the sequence codified by  $a$ by the formula: 

\[
\mathrm{Lgh}(a):=\left\{
\begin{array}{ccc}
0&\mathrm{if}&a\not\in\mathrm{Seq}\vee a=1\\
\mu x\leq a[p_{x}|a\wedge \neg(p_{x+1}|a)]&\mathrm{if}&a\in\mathrm{Seq}\wedge a\neq 1.
\end{array}\right.
\]

\noindent  `$(a)_{x}$' means that `$x$ is the $x$-th component of the sequence codified by $a$' and it is defined as:

\[
(a)_{x}:=\mu y\leq x+1[p_{x}^{y+1}|a\wedge \neg(p_{x}^{y+2}|a)].
\]
\end{defi}

\begin{remark}\rm Let be $a\in$Seq then we can write $a=\prod_{i\leq\mathrm{Lgh}(a)}p_{i}^{(a)_{i}+1}$. So, we define \textbf{concatenation} of sequences $a$ and $b$ by 

\[
a*b=a\prod_{x\leq\mathrm{Lgh}(b)}p_{\mathrm{Lgh}(a)+x+1}^{(b)_{x}+1},\hspace{0.2mm}\mathrm{if}\hspace{0.2mm}a, b\neq1.
\]
\end{remark}

Gödel introduced a pairing function $<, >$ which assigns to each $(x, y)\in\mathbb{N}\times\mathbb{N}$ a unique single natural number. This will be useful later for encoding partitions of finite subsets of natural numbers. Gödel pairing function is given by

\[
\begin{array}{cccc}
<, >:&\mathbb{N}\times\mathbb{N}&\longrightarrow&\mathbb{N}\\
 & (x, y)&\longmapsto&\frac{(x+y)(x+y+1)}{2}+y.
\end{array}
\]

\begin{lemma}
For any four given natural numbers $x, y, u, v$ the following is true $<x, y>=<u, v>$ iff $x=y$ and $u=v$. 
\end{lemma}

Let $m, n, c\in\mathbb{N}$ with $n\leq m$. We write `$H\in\mathrm{Part}([m]^{n}, c)$' to mean `$H$ is a partition of $[m]^{n}$ into $c$ pieces'. There exists a formula which it associates to  $H$ an unique natural number $\lceil H\rceil$.  In fact:

\[
\lceil H\rceil:=\prod_{\substack{1\leq r\leq n\\0\leq \alpha\leq c-1}}p_{r}^{\frac{\left[\displaystyle\prod_{\substack{1\leq_{i_{j}}\leq m\\0\leq j\leq n-1}}p_{i_{j}}^{m-i_{j}}+c_{\alpha}+1\right]\left[\displaystyle\prod_{\substack{1\leq_{i_{j}}\leq m\\0\leq j\leq n-1}}p_{i_{j}}^{m-i_{j}}+c_{\alpha}\right]}{2}+c_{\alpha} }.
\]

 We can write a subset of $[m]^{n}$ as $\{m-i_{0}, ..., m-i_{n-1}\}$ been $1\leq i_{j}\leq m$ for $0\leq j\leq n-1$. Now, $\mathrm{card}([m]^{n})=\binom{m}{n}=\frac{m!}{n!(m-n)!}$. Then we can list the elements of $[m]^{n}$ as: $[m]^{n}=m_{1}, ..., m_{\binom{m}{n}}$. Let be $H\in\mathrm{Part}([m]^{n}, c)$. We can write $H$ as  
 \\ $H:=\left\{(m_{\delta}, c_{\alpha}):1\leq\delta\leq\binom{m}{n}\wedge 0\leq\alpha\leq\ c\right\}$. Let $m_{\delta}\in[m]^{n}$, for $1\leq \delta\leq \binom{m}{n}$. Then

\begin{eqnarray*}
\lceil m_{\delta}\rceil&=&\lceil\{m-i_{0}, ..., m-i_{n-1}\}\rceil\\
&=&\prod_{\substack{1\leq_{i_{j}}\leq m\\0\leq j\leq n-1}}p_{i_{j}}^{m-i_{j}}
\end{eqnarray*}

\noindent for certain $1\leq i_{0}, ..., i_{n-1}\leq m$. So, for each pair $(m_{\delta}, c_{\alpha})$ we have a unique corresponding  pair 

\[
\left(\prod_{\substack{1\leq_{i_{j}}\leq m\\0\leq j\leq n-1}}p_{i_{j}}^{m-i_{j}}, c_{\alpha} \right)
\]

\noindent Then

\begin{eqnarray*}
\lceil m_{\delta}, c_{\alpha}\rceil&=&\left\langle\lceil m_{\delta}\rceil, c_{\alpha}\right\rangle\\
&=&\left\langle\lceil\{m-i_{0}, ..., m-i_{n-1}\}\rceil, c_{\alpha}\right\rangle\\
&=&\left\langle\prod_{\substack{1\leq_{i_{j}}\leq m\\0\leq j\leq n-1}}p_{i_{j}}^{m-i_{j}}, c_{\alpha} \right\rangle\\
&=&\frac{\left[\displaystyle\prod_{\substack{1\leq_{i_{j}}\leq m\\0\leq j\leq n-1}}p_{i_{j}}^{m-i_{j}}+c_{\alpha}+1\right]\left[\displaystyle\prod_{\substack{1\leq_{i_{j}}\leq m\\0\leq j\leq n-1}}p_{i_{j}}^{m-i_{j}}+c_{\alpha}\right]}{2}+c_{\alpha} 
\end{eqnarray*}

\noindent Finally

\begin{equation*}\label{partitions_code}
\lceil H\rceil: =\prod_{\substack{1\leq r\leq n\\0\leq \alpha\leq c-1}}p_{r}^{\frac{\left[\displaystyle\prod_{\substack{1\leq_{i_{j}}\leq m\\0\leq j\leq n-1}}p_{i_{j}}^{m-i_{j}}+c_{\alpha}+1\right]\left[\displaystyle\prod_{\substack{1\leq_{i_{j}}\leq m\\0\leq j\leq n-1}}p_{i_{j}}^{m-i_{j}}+c_{\alpha}\right]}{2}+c_{\alpha} }.
\end{equation*}

\noindent \textbf{*PH as a sentence on sequences.} 
\medskip

\noindent\textbf{PH.} For any natural numbers $b, e, y$ codifying sequences of length $n, c$ and at least  $\lambda$, respectively, there exists a natural number $a$ codifying a sequence of length $m$ such that,  for any $H\in\mathrm{Part}([m]^{n}, c)$, there exits a unique natural number $\beta\leq\mathrm{Lgh}(e)$ such that for any $1\leq k_{1}\leq\binom{m}{n}$, $(m_{k_{1}}, c_{\beta})\in H$ and there doesn't exists $\alpha\in\mathbb{N}$ such that for some $0\leq k_{2}\leq\binom{m}{n}$, $(m_{k_{2}}, c_{\alpha})\in H$  whenever $b$ is such that there is a sequence $s$ codifying a sequence such that $y=b*s$ and it holds that $\mathrm{Lgh}(y)\geq(y)_{0}$. In symbols:
\medskip

\hrule
\medskip

\noindent $\forall b\in\mathrm{Seq}\forall e\in\mathrm{Seq}\forall y\in\mathrm{Seq}\exists a\in\mathrm{ Seq}[[\mathrm{Lgh}(b)=n\wedge\mathrm{Lgh}(a)=m\wedge\mathrm{Lgh}(e)=c\wedge\mathrm{Lgh}(y)\leq\lambda\wedge n\leq m]:(\forall H\in\mathrm{Part}([m]^{n}, c))[(s\exists\in\mathrm{Seq})[y=b*s]\to(\exists!\beta\leq\mathrm{Lgh}(e)\in\mathbb{N}))[(\forall k_{1}\leq\binom{m}{n}\in\mathbb{N})[(m_{k_{1}}, c_{\beta})\in H]\wedge(\neg\exists\alpha\leq\mathrm{Lgh}(e)\in\mathbb{N})[(\exists k_{2} \leq \binom{m}{n}\in\mathbb{N})[(m_{k_{2}}, c_{\alpha})\in H]]]] \wedge\\ \mathrm{Lgh}(y) \geq (y)_{0}]$.
\medskip

\hrule
\medskip

This is now a formulation in $\mathcal{L}_{A}$ for PH. We can see is term used here has an explicit formulation in arithmetic. So, it we want to have pure symbolic formulation with symbols of $\mathcal{L}_{A}$ we would have to take the meaning of Seq, Lgh, etc. and put it into the formula above.  The formula we obtain after doing this is the following. You must note it can be even longer!

\newpage

\noindent $
 [(\forall b(b=1\vee(b>1\wedge\forall x_{1}\leq b[p_{x_{1}+1}|b\to p_{x_{1}}|b])))
  (\forall e(e=1\vee(e>1\wedge\forall x_{2}\leq e[p_{x_{2}+1}|e\to p_{x_{2}}|e])))
  (\forall y(y=1\vee(y>1\wedge\forall x_{3}\leq y[p_{x_{3}+1}|y\to p_{x_{3}}|y])))
  (\exists a(a=1\vee(a>1\wedge\forall x_{4}\leq a[p_{x_{4}+1}|a\to p_{x_{4}}|a])))]\\
 \left\{\begin{array}{c}
 \\
  \\
\end{array}\right
 [
 \mu x_{5}\leq b[p_{x_{5}}|b\wedge\neg(p_{x_{5}+1}|b)])=n\wedge\\
 \mu x_{6}\leq a[p_{x_{5}}|a\wedge\neg(p_{x_{5}+1}|a)])=m\wedge
 \mu x_{7}\leq e[p_{x_{5}}|e\wedge\neg(p_{x_{5}+1}|e)])=c\wedge\\
 \mu x_{8}\leq y[p_{x_{5}}|y\wedge\neg(p_{x_{5}+1}|y)])\leq\lambda\wedge n\leq m
 ]:\\
 \left[
 \forall H=\displaystyle\prod_{\substack{0\leq r\leq n-1\\0\leq \alpha\leq c-1}}p_{r}^{\frac{\left[\displaystyle\prod_{\substack{1\leq_{i_{j}}\leq m\\0\leq j\leq n-1}}p_{i_{j}}^{m-i_{j}}+c_{\alpha}+1\right]\left[\displaystyle\prod_{\substack{1\leq_{i_{j}}\leq m\\0\leq j\leq n-1}}p_{i_{j}}^{m-i_{j}}+c_{\alpha}\right]}{2}+c_{\alpha}}
 \right]\\
 (
 \exists s(s=1\vee(s>1\wedge\forall x_{6}\leq s[p_{x_{6}+1}|s\to p_{x_{6}}|s]))\\
\left(y=b\displaystyle\prod_{x_{7}\leq\mu x_{8}\leq s[p_{x_{8}}|s\wedge\neg(p_{x_{8}+1}|s)]}p_{\mu x_{5}\leq y[p_{x_{5}}|y\wedge\neg(p_{x_{5}+1}|y)]+x_{7}+1}^{\mu x_{9}\leq x_{7}+1[p_{x_{7}}^{x_{9}+1}|s\wedge\neg(p_{x_{7}}^{x_{9}+1}|s)]+1}\right)\to\\
(\exists!\beta=\underbrace{1+\cdots+1}_{\beta-\mathrm{times}}\leq\mu x_{7}\leq e[p_{x_{7}}|e\wedge\neg(p_{x_{7}+1}|e)]))\\
\left[\left(\forall k_{1}=\underbrace{1+\cdots+1}_{k_{1}-\mathrm{times}}\leq\\
\frac{m(m-1)(m-2)\cdots(m-(m-1))}{n(n-1)\cdots(n-(n-1))(m-n)(m-n-1)\cdots((m-n-(m-n-1)))}=h\right)\right.\\\
\exists\eta=\underbrace{1+\cdots+1}_{\eta-\mathrm{times}}\leq h
\left[p^{\lceil \lceil m_{k_{1}}\rceil, c_{\beta}\rceil }_{\eta}|H \right]\wedge\\(\neg\exists \alpha=\underbrace{1+\cdots+1}_{a-\mathrm{times}}\leq\mu x_{7}\leq e[p_{x_{5}}|e\wedge\neg(p_{x_{5}+1}|e)]))\\
\left[\left(\exists k_{2}=\underbrace{1+\cdots+1}_{k_{2}-\mathrm{times}}\leq\frac{m(m-1)(m-2)\cdots(m-(m-1))}{n(n-1)\cdots(n-(n-1))(m-n)(m-n-1)\cdots((m-n-(m-n-1)))}\right)\right.\\
\exists\rho=\underbrace{1+\cdots+1}_{\rho-\mathrm{times}}\leq h
\left.\left.\left.\left(p_{\rho}^{\lceil\lceil m_{k_{2}}\rceil, c_{\alpha}\rceil}|H\right)\right]\right]\right]\wedge\\
\mu x_{8}\leq y[p_{x_{5}}|y\wedge\neg(p_{x_{5}+1}|y)])\geq\mu x_{10}\leq 0+1[p_{x_{10}}^{0+1}|y\wedge \neg(p_{x_{10}}^{0+2}|y)]
\left.\begin{array}{c}
 \\
  \\
\end{array}\right\}
$

For getting this formulation, we wrote $\lceil\lceil m_{k_{2}}\rceil, c_{\alpha}\rceil$ to mean the code of the pair $(m_{k_{2}}, c_{\alpha})$. We also used the notation $p_{n}$ to mean $p_{n}$ is the $n$-th prime number. We also omitted replace the meaning of the $\mu$-bounded operator. But, note that with doing this we don't lose the sense of the formula, but we agree that it is easier to read.

\end{document}